# Spectrum of The Direct Sum of Operators

by

**E.OTKUN ÇEVİK and Z.I.ISMAILOV**

Karadeniz Technical University, Faculty of Sciences, Department of Mathematics
61080 Trabzon, TURKEY
e-mail adress : zameddin@yahoo.com

**Abstract:** In this work, a connection between some spectral properties of direct sum of operators in the direct sum of Hilbert spaces and its coordinate operators has been investigated.

**Keywords:** Direct sum of Hilbert spaces and operators; spectrum and resolvent sets; continuous and compact operators; discreteness of spectrum; asymptotics of eigenvalues;

**2000 AMS Subject Classification:** 47A10.

## 1. Introduction

It is known that infinite direct sum of Hilbert spaces $H_n$, $n \geq 1$ and infinite direct sum of operators $A_n$ in $H_n$, $n \geq 1$ are defined as

$$H = \bigoplus_{n=1}^{\infty} H_n = \left\{ u = (u_n) \: : \: u_n \in H_n \, , \, n \geq 1 \, , \, \|u\|_H^2 = \sum_{n=1}^{\infty} \|u_n\|_{H_n}^2 < \infty \right\}$$

and

$$A = \bigoplus_{n=1}^{\infty} A_n \, , \, D(A) = \left\{ u = (u_n) \in H \: : \: u_n \in D(A_n) \, , \, n \geq 1, \, Au = (A_n u_n) \in H \right\} ,$$

$$A : D(A) \subset H \to H \quad \text{(see [1])}.$$

The general theory of linear closed operators in Hilbert spaces and its applications to physical problems has been investigated by many mathematicians (for example, see [1]).

However, many physical problems of today arising in the modelling of processes of multiparticle quantum mechanics, quantum field theory and in the physics of rigid bodies support to study a theory of linear direct sum of operators in the direct sum of Hilbert spaces (see [2]-[6] and references in it).

In this paper, a connection between spectrum, resolvent sets, discreteness of the spectrum (sec. 2) and asymtotical behaviour of the eigenvalues (sec. 3) of direct sum of operators defined in the direct sum of Hilbert spaces and suitable properties of coordinate operators has been established. The obtained results has been supported with applications.

These and related problems in the case continuous direct sum of the Hilbert space operators have been investigated in works [7 − 10]. But in these works has not been considered a connection between parts of the spectrum of direct sum operator and suitable parts of the spectrum their coordinate operators. In this paper given sharp formulaes in the this sense.

## 2. On the spectrum of direct sum of operators

In this section, the relationship between the spectrum and resolvent sets of the direct sum of operators and its coordinate operators will be investigated.

First of all it will be investigated the continuity and compactness properties of the operator $A = \bigoplus_{n=1}^{\infty} A_n$ in $H = \bigoplus_{n=1}^{\infty} H_n$ in case when $A_n \in L(H_n)$ for each $n \geq 1$.

It is easy to see that the following propositions are true in general.

**Theorem 2.1.** Let $A = \bigoplus_{n=1}^{\infty} A_n$, $H = \bigoplus_{n=1}^{\infty} H_n$ and for any $n \geq 1$ $A_n \in L(H_n)$. In order for $A \in L(H)$ the necessary and sufficient condition is $\sup_{n \geq 1} \|A_n\| < \infty$.

In addition, in this case when $A \in L(H)$ it is true $\|A\| = \sup_{n \geq 1} \|A_n\|$ (see [7]).

**Theorem 2.2.** Let $A_n \in C_{\infty}(H_n)$ for each $n \geq 1$. In this case $A = \bigoplus_{n=1}^{\infty} A_n \in C_{\infty}(H)$ if and only if $\lim_{n \to \infty} \|A_n\| = 0$.

Furthermore, the following main result can be proved.

**Theorem 2.3.** For the parts of spectrum and resolvent sets of the operator $A = \bigoplus_{n=1}^{\infty} A_n$ in Hilbert space $H = \bigoplus_{n=1}^{\infty} H_n$ the following claims are true

$$\sigma_p(A) = \bigcup_{n=1}^{\infty} \sigma_p(A_n),$$

$$\sigma_c(A) = \left\{ \left( \bigcup_{n=1}^{\infty} \sigma_p(A_n) \right)^c \cap \left( \bigcup_{n=1}^{\infty} \sigma_r(A_n) \right)^c \cap \left( \bigcup_{n=1}^{\infty} \sigma_c(A_n) \right) \right\}$$

$$\bigcup \left\{ \lambda \in \bigcap_{n=1}^{\infty} \rho(A_n) : sup\|R_\lambda(A_n)\| = \infty \right\},$$

$$\sigma_r(A) = \left(\bigcup_{n=1}^{\infty} \sigma_p(A_n)\right)^c \cap \left(\bigcup_{n=1}^{\infty} \sigma_r(A_n)\right),$$

$$\rho(A) = \{\lambda \in \bigcap_{n=1}^{\infty} \rho(A_n): sup\|R_\lambda(A_n)\| < \infty\}$$

**Proof.** The validity of first claim of given relations is clear. Moreover, it is easy to prove the fourth equality using the theorem 2.1.

Now we prove the second relation on the continuous spectrum.

Let $\lambda \in \sigma_c(A)$. In this case by the definition of continuous spectrum $A - \lambda E$ is a one-to-one operator, $R(A - \lambda E) \neq H$ and $\overline{R(A - \lambda E)}$ is dense in $H$. Consequently, for any $n \geq 1$ an operator $A_n - \lambda E_n$ is a one-to-one operator in $H_n$, there exists $m \in \square$ such that $R(A_m - \lambda E_m) \neq H_m$ and for any $n \geq 1$ linear manifold $R(A_n - \lambda E_n)$ is dense in $H_n$ or $\lambda \in \rho(A_m)$ for each $m \geq 1$ but $sup\{\|R_\lambda(A_m)\|: m \geq 1\} = \infty$. This means that

$$\lambda \in \left\{\left(\bigcap_{n=1}^{\infty}[\sigma_c(A_n) \cup \rho(A_n)]\right) \cap \left(\bigcup_{n=1}^{\infty} \sigma_c(A_n)\right)\right\}$$

$$\bigcup \left\{\lambda \in \bigcap_{n=1}^{\infty} \rho(A_n): sup\|R_\lambda(A_n)\| = \infty\right\}$$

On the contrary, now suppose that for the point $\lambda \in \square$ the above relation is satisfied. Consequently, either for any $n \geq 1$

$$\lambda \in \sigma_c(A_n) \cup \rho(A_n),$$

or $\lambda \in \bigcap_{n=1}^{\infty} \rho(A_n): sup\|R_\lambda(A_n)\| = \infty$,

and there exist $m \in \square$ such that

$$\lambda \in \sigma_c(A_m).$$

That is, for any $n \geq 1$ $A_n$ is a one-to-one operator, $\overline{R(A_n - \lambda E_n)} = H_n$ and $R(A_m - \lambda E_m) \neq H_m$. And from this it implies that the operator $A = \bigoplus_{n=1}^{\infty} A_n$ is a one-to-one operator, $R(A - \lambda E) \neq H$ and $\overline{R(A - \lambda E)} = H$. Hence $\lambda \in \sigma_c(A_m)$.

On the other hand the simple calculations show that

$$\left[\bigcap_{n=1}^{\infty}(\sigma_c(A_n)\cup\rho(A_n))\right]\cap\left[\bigcup_{n=1}^{\infty}\sigma_c(A_n)\right]=\left[\bigcup_{n=1}^{\infty}\sigma_p(A_n)\right]^c\cap\left[\bigcup_{n=1}^{\infty}\sigma_r(A_n)\right]^c\cap\left(\bigcup_{n=1}^{\infty}\sigma_c(A_n)\right).$$

By the similarly idea can be proved the validity of the third equality of the theorem.

**Example 2.4.** Consider the following multipoint differential operator for first order

$$A_n u_n = u_n'(t) \ , \ H_n = L^2(\Delta_n) \ , \ \Delta_n = (a_n, b_n) \ , \ -\infty < a_n < b_n < a_{n+1} < \ldots < +\infty$$

$$A_n : D(A_n) \subset H_n \to H_n \ , \ D(A_n) = \{u_n \in W_2^1(\Delta_n) : u_n(a_n) = u_n(b_n)\} \ , \ n \geq 1 \ ,$$

$$A = \bigoplus_{n=1}^{\infty} A_n \quad \text{and} \quad H = \bigoplus_{n=1}^{\infty} H_n \ .$$

For any $n \geq 1$ operator $A_n$ and $A$ are normal, $\sigma(A_n) = \sigma_p(A_n) = \left\{\dfrac{2k\pi i}{b_n - a_n} : k \in \mathbb{Z}\right\}$ and eigenvectors according to the eigenvalue $\lambda_{nk}$, $n \geq 1$, $k \in \mathbb{Z}$ are in the form

$$u_{nk}(t) = c_{nk} \exp(\lambda_{nk}(t - a_n)) \ , \ t \in \Delta_n \ , \ c_{nk} \in \mathbb{C} \setminus \{0\} \ [11].$$

In this case $\displaystyle\sum_{n=1}^{\infty}\|u_{nk}\|_{H_n}^2 = \sum_{n=1}^{\infty}\int_{\Delta_n}|c_{nk}|^2 \left|\exp(\lambda_{nk}(t-a_n))\right|^2 dt = \sum_{n=1}^{\infty}|c_{nk}|^2(b_n - a_n).$

The coefficients $c_{nk}$ may be choosen such that the last series to be convergent. This means that $\lambda_{nk} \in \sigma_p(A)$. From this and Theorem 2.1 it is obtained that $\sigma_p(A) = \displaystyle\bigcup_{n=1}^{\infty}\sigma_p(A_n)$.

**Definition 2.5[12].** Let $T$ be a linear closed and densely defined operator in any Hilbert space $H$. If $\rho(T) \neq \varnothing$ and for $\lambda \in \rho(T)$ the resolvent operator $R_\lambda(T) \in C_\infty(H)$, then operator $T : D(T) \subset H \to H$ is called a operator with discrete spectrum.

Note that the following results are true.

It is clear that if the operator $A = \displaystyle\bigoplus_{n=1}^{\infty} A_n$ is an operator with discrete spectrum in $H = \displaystyle\bigoplus_{n=1}^{\infty} H_n$, then for every $n \geq 1$ the operator $A_n$ is also in $H_n$.

The following proposition is proved by using the theorem 2.2.

**Theorem 2.6.** If $A = \bigoplus_{n=1}^{\infty} A_n$, $A_n$ is an operator with discrete spectrum in $H_n$, $n \geq 1$, $\bigcap_{n=1}^{\infty} \rho(A_n) \neq \varnothing$ and $\lim_{n \to \infty} \|R_\lambda(A_n)\| = 0$, then $A$ is an operator with discrete spectrum in $H$.

**Proof:** In this case for each $\lambda \in \bigcap_{n=1}^{\infty} \rho(A_n)$ we have $R_\lambda(A_n) \in C_\infty(H_n)$, $n \geq 1$.

Now we define the operator $K := \bigoplus_{n=1}^{\infty} R_\lambda(A_n)$ in $H$. In this case for every $u = (u_n) \in D(A)$ we have

$$K(A - \lambda E)(u_n) = \bigoplus_{n=1}^{\infty} R_\lambda(A_n) \left( \bigoplus_{n=1}^{\infty} (A_n - \lambda E_n) \right)(u_n) = \bigoplus_{n=1}^{\infty} R_\lambda(A_n)\left((A_n - \lambda E_n)u_n\right)$$
$$= \left(R_\lambda(A_n)(A_n - \lambda E_n)u_n\right) = (u_n)$$

and

$$(A - \lambda E)K(u_n) = (A - \lambda E)\left( \bigoplus_{n=1}^{\infty} R_\lambda(A_n) \right)(u_n) = (A - \lambda E)\left(R_\lambda(A_n)u_n\right)$$
$$= \left( \bigoplus_{n=1}^{\infty} (A_n - \lambda E_n) \right)\left(R_\lambda(A_n)u_n\right) = \left((A_n - \lambda E_n)R_\lambda(A_n)u_n\right) = (u_n)$$

These relations show that $R_\lambda(A) = \bigoplus_{n=1}^{\infty} R_\lambda(A_n)$

Furthermore, we define the following operators $K_m : H \to H$, $m \geq 1$ in the form

$$K_m u := \{R_\lambda(A_1)u_1, R_\lambda(A_2)u_2, \ldots, R_\lambda(A_m)u_m, 0, 0, \ldots\}, u = (u_n) \in H.$$

Now the convergence in operator norm of the operators $K_m$ to the operator $K$ will be investigated. For the $u = (u_n) \in H$ we have

$$\|K_m u - Ku\|_H^2 = \sum_{n=m+1}^{\infty} \|R_\lambda(A_n)u_n\|_{H_n}^2 \leq \sum_{n=m+1}^{\infty} \|R_\lambda(A_n)\|^2 \|u_n\|_{H_n}^2 \leq \left( \sup_{n \geq m+1} \|R_\lambda(A_n)\| \right)^2 \sum_{n=1}^{\infty} \|u_n\|_{H_n}^2$$
$$= \left( \sup_{n \geq m+1} \|R_\lambda(A_n)\| \right)^2 \|u\|_H^2$$

From this it is obtained that $\|K_m - K\| \leq \sup_{n \geq m+1} \|R_\lambda(A_n)\|$, $m \geq 1$.

This means that sequence of the operators $(K_m)$ converges in operator norm to the operator $K$. Then by the important theorem of the theory of compact operators $K \in C_\infty(H)$ [1], because for any $m \geq 1$ $K_m \in C_\infty(H)$.

**Example 2.7.** Consider that the following family of the operators in the form

$$A_n := \frac{d}{dt} + S_n, \quad S_n^* = S_n > 0, \quad S_n^{-1} \in C_\infty(H),$$

$$A_n : D(A_n) \subset L_n^2 \to L_n^2, \quad \Delta_n = (a_n, b_n), \quad \sup_{n \geq 1}(b_n - a_n) < \infty,$$

$$D(A_n) = \{u_n \in W_2^1(H, \Delta_n) : u_n(b_n) = W_n u_n(a_n), \; A_n^{-1} W_n = W_n A_n^{-1}\},$$

where $L_n^2 = L^2(H, \Delta_n)$, $n \geq 1$, $H$ is any Hilbert space and $W_n$ is a unitary operator in $H$, $n \geq 1$ (for this see [11]). For any $n \geq 1$ an operator $A_n$ is normal with discrete spectrum and $\bigcap_{n=1}^{\infty} \rho(A_n) \neq \varnothing$.

For the $\lambda \in \bigcap_{n=1}^{\infty} \rho(A_n)$ and sufficiently large $n \geq 1$ a simple calculation shows that

$$R_\lambda(A_n) f_n(t) = e^{-(S_n - \lambda E_n)(t - a_n)} \left( E - W_n^* e^{-(S_n - \lambda E_n)(b_n - a_n)} \right)^{-1} W_n^* \int_{\Delta_n} e^{-(S_n - \lambda E_n)(b_n - s)} f_n(s) ds$$

$$+ \int_{a_n}^{t} e^{-(S_n - \lambda E_n)(t - s)} f_n(s) ds, \quad f_n \in L_n^2, \; n \geq 1$$

On the other hand the following estimates are true

$$\left\| \int_{a_n}^{t} e^{-(S_n - \lambda E_n)(t - s)} f_n(s) ds \right\|_{L_n^2}^2 \leq \int_{\Delta_n} \left( \int_{a_n}^{t} \left\| e^{-(S_n - \lambda E_n)(t - s)} \right\| \|f_n(s)\|_H ds \right)^2 dt$$

$$\leq \int_{\Delta_n} \left( \int_{a_n}^{t} \left\| e^{-(S_n - \lambda E_n)(t - s)} \right\|^2 ds \right) dt \int_{\Delta_n} \|f_n(s)\|_H^2 ds = \int_{\Delta_n} \left( \int_{a_n}^{t} \left\| e^{-(S_n - \lambda_r E_n)(t - s)} \right\|^2 ds \right) dt \|f_n\|_{L_n^2}^2$$

$$= \int_{\Delta_n} \left( \int_{a_n}^{t} e^{2\lambda_r(t - s)} \left\| e^{-S_n(t - s)} \right\|^2 ds \right) dt \|f_n\|_{L_n^2}^2 = \int_{\Delta_n} \left( \int_{a_n}^{t} e^{2(\lambda_r - \lambda_1^{(n)})(t - s)} ds \right) dt \|f_n\|_{L_n^2}^2$$

$$= \frac{1}{4(\lambda_r - \lambda_1^{(n)})^2} \left[ 2(\lambda_r - \lambda_1^{(n)})(a_n - b_n) - 1 + e^{2(\lambda_r - \lambda_1^{(n)})(b_n - a_n)} \right] \|f_n\|_{L_n^2}^2, \qquad (2.1)$$

$$\left\|\left(E-W_n^*e^{-(S_n-\lambda E_n)(b_n-a_n)}\right)^{-1}W_n^*\right\| = \left\|\sum_{m=0}^{\infty}\left(W_n^*e^{-(S_n-\lambda E_n)(b_n-a_n)}\right)^m\right\| \leq \sum_{m=0}^{\infty}\left\|e^{-(S_n-\lambda_r E_n)(b_n-a_n)}\right\|^m$$

$$=\left(1-\left\|e^{-(S_n-\lambda_r E_n)(b_n-a_n)}\right\|\right)^{-1} = \left(1-e^{(\lambda_r-\lambda_1^{(n)})(b_n-a_n)}\right)^{-1}, \tag{2.2}$$

$$\left\|\int_{\Delta_n} e^{-(S_n-\lambda E_n)(b_n-s)} f_n(s)\,ds\right\|^2 \leq \frac{1}{2(\lambda_r-\lambda_1^{(n)})}\left[e^{2(\lambda_r-\lambda_1^{(n)})(b_n-a_n)}-1\right]\|f_n\|_{L_n^2}^2, \tag{2.3}$$

Hence from $(2.2)$ and $(2.3)$ we have

$$\left\|e^{-(S_n-\lambda E_n)(t-a_n)}\left(E-W_n^*e^{-(S_n-\lambda E_n)(b_n-a_n)}\right)^{-1}W_n^*\int_{\Delta_n} e^{-(S_n-\lambda E_n)(b_n-s)} f_n(s)\,ds\right\|^2_{L_n^2} \leq$$

$$\leq \int_{\Delta_n} e^{2\lambda_r(t-a_n)}\left\|e^{-S_n(t-a_n)}\right\|^2 dt \left\|\left(E-W_n^*e^{-(S_n-\lambda E_n)(b_n-a_n)}\right)^{-1}W_n^*\right\|^2 \left\|\int_{\Delta_n} e^{-(S_n-\lambda E_n)(b_n-s)} f_n(s)\,ds\right\|^2_{L_n^2} \leq$$

$$\leq \frac{1}{4\lambda_r(\lambda_r-\lambda_1^{(n)})}\left(e^{2\lambda_r(b_n-a_n)}-1\right)\left(1-e^{(\lambda_r-\lambda_1^{(n)})(b_n-a_n)}\right)^{-1}\left(e^{2(\lambda_r-\lambda_1^{(n)})(b_n-a_n)}-1\right)\|f_n\|_{L_n^2}^2, \tag{2.4}$$

where, $\lambda_r$ is the real part of $\lambda$ and $\lambda_1^{(n)}$ is the first eigenvalue of the operator $S_n$, $n\geq 1$.

Therefore, from estimates $(2.1)$ and $(2.4)$ the following result is obtained.

**Proposition 2.8.** If $\lambda \in \bigcap_{n=1}^{\infty}\rho(A_n)$, $\sup_{n\geq 1}(b_n-a_n)<\infty$ and $\lambda_1^{(n)}(S_n)\to\infty$ as $n\to\infty$, then

$$\|R_\lambda(A_n)\|\to 0 \text{ as } n\to\infty.$$

Consequently, the operator $A = \bigoplus_{n=1}^{\infty} A_n$ is an operator with discrete spectrum in $L^2 = \bigoplus_{n=1}^{\infty} L_n^2$.

## 3. Asymptotical behaviour of the eigenvalues

In this section asymptotical behaviour for the eigenvalues of the operator $A = \bigoplus_{n=1}^{\infty} A_n$ in $H = \bigoplus_{n=1}^{\infty} H_n$ will be investigated in a special case.

**Theorem 3.1.** Assume that teh operator $A$ in $H$ and $A_n$ in $H_n$, $n \geq 1$ are operators with discrete spectrum and for $i, j \geq 1$, $i \neq j$, $\sigma(A_i) \cap \sigma(A_j) = \emptyset$ is satisfied. If $\lambda_m(A_n) \Box c_n m^{\alpha_n}$, $0 < c_n, \alpha_n < \infty$, $m \to \infty$, $\sum_{n=1}^{\infty} c_n^{-1/\alpha_n} < \infty$, $n \geq 1$ and there exists $q \in \Box$ such that $\alpha_q = \inf_{n \geq 1} \alpha_n > 0$, then $\lambda_n(A) \Box \gamma n^{\alpha}$, $0 < \gamma, \alpha = \alpha_q < \infty$ as $n \to \infty$.

**Proof:** First of all note that by the Theorem 2.3 $\sigma_p(A) = \bigcup_{n=1}^{\infty} \sigma_p(A_n)$.

Here it is denoted by $N(T; \lambda) := \sum_{|\lambda(T)| \leq \lambda} 1$, $\lambda \geq 0$, that is, a number of eigenvalues of the some linear closed operator $T$ in any Hilbert space with modules of these eigenvalues less than or equal to $\lambda$, $\lambda \geq 0$. This function takes values in the set of non-negative integer numbers and in case of unbounded operator $T$ it is nondecreasing and tends to $\infty$ as $\lambda \to \infty$.

Since for every $i, j \geq 1$, $i \neq j$, $\sigma(A_i) \cap \sigma(A_j) = \emptyset$, then $N(A; \lambda) = \sum_{n=1}^{\infty} N(A_n; \lambda)$

In this case it is clear that

$$\frac{N(A; \lambda)}{\lambda^{1/\alpha}} \Box \sum_{n=1}^{\infty} c_n^{-1/\alpha_n} \lambda^{\frac{1}{\alpha_n} - \frac{1}{\alpha}} = \sum_{n=1}^{\infty} c_n^{-1/\alpha_n} \lambda^{\frac{\alpha - \alpha_n}{\alpha \alpha_n}} = \sum_{n=1}^{\infty} c_n^{-1/\alpha_n} \left(\frac{1}{\lambda}\right)^{\frac{\alpha_n - \alpha}{\alpha \alpha_n}}, \quad \lambda \geq 0.$$

The last series is uniformly convergent in $(1, \infty)$ on $\lambda$. Then $\lim_{\lambda \to \infty} \sum_{n=1}^{\infty} c_n^{-1/\alpha_n} \left(\frac{1}{\lambda}\right)^{\frac{\alpha_n - \alpha}{\alpha \alpha_n}} = c_q^{-1/\alpha_q}$

Therefore $N(A; \lambda) \Box c \lambda^{1/\alpha}$, $0 < c = c_q^{-1/\alpha_q}, \alpha < \infty$ as $\lambda \to \infty$.

Then the following asymptotic behaviour of eigenvalues of the operator $A$ in $H$ $\lambda_n(A) \Box \gamma n^{\alpha}$, $0 < \gamma, \alpha < \infty$ as $n \to \infty$ is true.

**Remark 3.2.** If in the above theorem the coefficients $\alpha_n$, $n \geq 1$ satisfy the following condition $\inf_{n \geq 1} \alpha_n > 0$, then for every $0 < \alpha < \inf_{n \geq 1} \alpha_n$, $N(A; \lambda) = o(\lambda^{1/\alpha})$ as $\lambda \to \infty$

**Remark 3.3.** If the every finitely many sets of the family $\sigma(A_n)$, $n \geq 1$ in complex plane intersect in the finitely many points, then it can be proved that claim of the Theorem 3.1 is valid in this case too.

**Example 3.4:** Let $H = \bigoplus_{n=1}^{\infty} H_n$, $H_n = \ell^2(\Box)$, $n \geq 1$, $A_n : D(A_n) \subset H_n \to H_n$, $A_n(u_m) := (c_{nm} u_m)$, $u = (u_m) \in D(A_n)$, $c_{nm} \in \Box$, $c_{nm} \neq c_{km}$, $n \neq k$, $n, k, m \geq 1$, $c_{nm} \Box k_n m^{\alpha_n}$, $0 < k_n < \infty$,

$1 \leq \alpha_n < \infty$ as $m \to \infty$, $\sum_{n=1}^{\infty} k_n^{-1/\alpha_n}$ is convergent and there exists $q \in \mathbb{N}$ such that $\alpha_q = \inf_{n \geq 1} \alpha_n$.

In this case, for any $n \geq 1$ $A_n$ is a linear normal operator and $\sigma(A_n) = \sigma_p(A_n) = \overline{\bigcup_{m=1}^{\infty} \{c_{nm}\}}$

Now we obtain the resolvent operator of $A_n$. Let $\lambda \in \bigcap_{n=1}^{\infty} \rho(A_n)$. Then from the relation $(A_n - \lambda E_n)(u_m) = (v_m)$, $n \geq 1$, $(v_m) \in H_n$, i.e $c_{nm} u_m - \lambda u_m = v_m$, $m \geq 1$.

It is established that $u_m = \dfrac{v_m}{c_{nm} - \lambda}$, $m \geq 1$, i.e $R_\lambda(A_n)(v_m) = \left(\dfrac{v_m}{c_{nm} - \lambda}\right)$, $n \geq 1$.

On the other hand since $c_{nm} \square k_n m^{\alpha_n}$, $\alpha_n \geq 1$ as $m \to \infty$, then for any $v = (v_m) \in H_n$ we have

$$\|R_\lambda(A_n)(v_m)\|_{H_n}^2 = \sum_{m=1}^{\infty} \left|\dfrac{v_m}{c_{nm} - \lambda}\right|^2 \leq \sum_{m=1}^{\infty} \left|\dfrac{1}{c_{nm} - \lambda}\right|^2 \cdot \sum_{m=1}^{\infty} |v_m|^2 = \sum_{m=1}^{\infty} \left|\dfrac{1}{c_{nm} - \lambda}\right|^2 \cdot \|v\|_{H_n}^2$$

Consequently, for any $n \geq 1$

$$\|R_\lambda(A_n)\| \leq \left(\sum_{m=1}^{\infty} \left|\dfrac{1}{c_{nm} - \lambda}\right|^2\right)^{1/2} \tag{3.1}$$

Moreover, it is known that a resolvent operator $R_\lambda(A_n)$, $n \geq 1$ is compact if and only if $\dfrac{1}{c_{nm} - \lambda} \xrightarrow[m \to \infty]{} 0$ [13]. Since $\lambda \neq c_{nm}$, $n, m \geq 1$ and conditions on $c_{nm}$, then the last condition is satisfied. Hence for any $n \geq 1$ $R_\lambda(A_n) \in C_\infty(H_n)$.

On the other hand since the series $\sum_{n=1}^{\infty} k_n^{-1/\alpha_n}$ is convergent, then from the inequality (3.1) for the $R_\lambda(A_n)$, $n \geq 1$ it is easy to see that $\lim_{n \to \infty} \|R_\lambda(A_n)\| = 0$, $\lambda \in \bigcap_{n=1}^{\infty} \rho(A_n)$

Hence by the Theorem 2.6 for the $\lambda \in \bigcap_{n=1}^{\infty} \rho(A_n)$ it is established that $R_\lambda(A) \in C_\infty(H)$. Then by the Theorem 2.3 it is true that $\sigma_p(A) = \bigcup_{n=1}^{\infty} \sigma_p(A_n)$.

Furthermore, the validity of the relation $\sigma(A_i) \cap \sigma(A_j) = \emptyset$, $i, j \geq 1$, $i \neq j$ is clear.

Therefore by the Theorem 3.1 $\lambda_n(A) \square \gamma n^\alpha$, $0 < \gamma, \alpha < \infty$ as $n \to \infty$